\documentclass[oneside,notitlepage,12pt]{article}

\usepackage{amssymb}
\usepackage[leqno]{amsmath}
\usepackage{amsfonts}
\usepackage{amsopn}
\usepackage{amstext}
\usepackage{amsthm}

\usepackage[all]{xy}
\newdir{ >}{{}*!/-9pt/@{>}}
\usepackage[colorlinks, backref]{hyperref}

\usepackage{calrsfs}
\usepackage{fourier-orns}
\usepackage{hieroglf} 
\usepackage{clock} 

\providecommand{\cal}{\mathcal}
\renewcommand{\Bbb}{\mathbb}

\newenvironment{pf}{\begin{proof}}{\end{proof}}

\newcommand{\Ef}{{\cal{F}}}

\newcommand{\Qyu}{{\Bbb{Q}}}
\newcommand{\Err}{{\Bbb{R}}}

\newcommand{\lam}{{\lambda}}

\newcommand{\eps}{\varepsilon}
\renewcommand{\phi}{\varphi}
\renewcommand{\rho}{\varrho}

\newcommand{\rest}{\restriction}

\newcommand{\ntr}{{n\in\omega}}
\newcommand{\Ntr}{n\in{\Bbb{N}}}
\newcommand{\loe}{\leq}
\newcommand{\goe}{\geq}

\newcommand{\subs}{\subseteq}
\newcommand{\sups}{\supseteq}

\newcommand{\oraz}{\qquad\text{and}\qquad}

\newcommand{\by}{/}

\newtheorem{tw}{Theorem}[section]

\newtheorem{lm}[tw]{Lemma}
\newtheorem{prop}[tw]{Proposition}
\newtheorem{claim}[tw]{Claim}
\theoremstyle{definition}

\theoremstyle{remark}
\newtheorem{uwgi}[tw]{Remark}

\newcommand{\setof}[2]{\{#1\colon #2\}}
\newcommand{\Bigsetof}[2]{\Bigl\{#1\colon #2\Bigr\}}

\newcommand{\sett}[2]{\{#1\}_{#2}}
\newcommand{\sn}[1]{\{#1\}} 
\newcommand{\pair}[2]{(#1, #2)} 
\newcommand{\triple}[3]{\langle #1, #2, #3 \rangle} 
\newcommand{\map}[3]{#1\colon #2 \to #3} 
\newcommand{\img}[2]{#1[#2]} 

\providecommand{\nat}{\omega}

\newcommand{\ciag}[1]{{\sett{{#1}_n}{\ntr}}}

\newcommand{\anorm}{\|\cdot\|}
\newcommand{\norm}[1]{\|#1\|}

\newcommand{\ubal}[1]{\operatorname{B}_{#1}}

\newcommand{\fK}{{\mathfrak{K}}}

\newcommand{\Ban}{{\mathfrak B_1}}

\newcommand{\cmp}{\circ} 

\newcommand{\G}{\mathbb G}

\newcommand{\Gurarii}{Gurari\u\i}
\newcommand{\gurdef}{\mathfrak G} 
\newcommand{\uop}{\mathbf\Omega} 
\newcommand{\almap}[3]{#1 \colon #2 \rightsquigarrow #3}

\title{A universal operator on the \Gurarii\ space}
\author{
{\sc Joanna Garbuli\'nska-W\c egrzyn}\thanks{Research of the first author supported by the ESF Human Capital Operational
Programme grant 6/1/8.2.1./POKL/ 2009.}\\
{\small Institute of Mathematics,}
{\small Jan Kochanowski University (POLAND)}\\
{\small Faculty of Mathematics and Computer Science,} 
{\small Jagiellonian University (POLAND)}\\
{\small\texttt{jgarbulinska@ujk.edu.pl}}
\and
{\sc Wies{\l}aw Kubi\'s}\thanks{Research of the second author supported by the GA\v CR grant P 201/12/0290.}\\
{\small Institute of Mathematics, Academy of Sciences of the Czech Republic}\\
{\small Institute of Mathematics,}
{\small Jan Kochanowski University (POLAND)}\\
{\small\texttt{kubis@math.cas.cz}}
}

\begin{document}
\maketitle

\begin{abstract}
We construct a nonexpansive linear operator on the \Gurarii\ space that ``captures" all nonexpansive linear operators between separable Banach spaces.
Some additional properties involving its restrictions to finite-dimensional subspaces describe this operator uniquely up to an isometry.

\noindent
{\bf MSC (2010):} 
47A05, 
47A65, 
46B04. 

\noindent
{\bf Keywords:}
Isometrically universal operator, \Gurarii\ space, almost isometry.
\end{abstract}

\tableofcontents

\section{Introduction}

There exist at least two different notions of universal operators between Banach spaces (by an \emph{operator} we mean a bounded linear operator).
Perhaps the most popular one, due to Caradus~\cite{Caradus} is the following: An operator $\map U X X$ is \emph{universal} if for every other operator $\map T X X$ there exist a $U$-invariant subspace $Y \subs X$ and a linear isomorphism $\map \phi X Y$ such that $\lam T = \phi^{-1} \cmp (U\rest Y) \cmp \phi$ for some constant $\lam > 0$.
Caradus~\cite{Caradus} described universal operators on the separable Hilbert space.
Some arguments from dilation theory show that the left-shift on the Hilbert space is actually universal in a stronger sense: the isomorphism $\phi$ is a linear isometry, whenever the operator $T$ is contractive and satisfies $\lim_{n\to\infty}T^n x = 0$ for every $x \in H$.
The details can be found in~\cite{AmbMul}.
A much weaker notion of a universal operator is due to Lindenstrauss and Pe\l czy\'nski~\cite{LinPel}: An operator $U$ is \emph{universal} for a given class $\Ef$ of operators if for every $T \in \Ef$ there exist operators $L, R$ such that $L \cmp T \cmp R = U$.
One of the results of \cite{LinPel} says that the ``partial sums'' operator $\map U {\ell_1}{\ell_\infty}$ is universal for the class of non-compact operators.

We are concerned with a natural concept of an \emph{isometrically universal} operator, that is, an operator $U$ between separable Banach spaces having the property that for every other operator acting on separable Banach spaces, whose norm does not exceed the norm of $U$, there exist isometric embeddings $i$, $j$ such that $U \cmp i = j \cmp T$.
This property is weaker than the isometric variant of Caradus' concept (since we allow two different embeddings and no invariant subspace), although much stronger than the universality in the sense of Lindenstrauss and Pe\l czy\'nski.

Our main result is the existence of an isometrically universal operator $\uop$.
We also formulate its extension property which describes this operator uniquely, up to isometries.
This is in contrast with the result of Caradus, where a rather general criterion for being universal is given.
It turns out that both the domain and the co-domain of our operator $\uop$ are isometric to the \Gurarii\ space.

Recall that the \emph{\Gurarii\ space} is the unique separable Banach space $\G$ satisfying the following condition: Given finite-dimensional Banach spaces $X \subs Y$, $\eps > 0$, every isometric embedding $\map i X \G$ extends to an $\eps$-isometric embedding $\map j Y \G$.
This space was constructed by \Gurarii~\cite{gurarii} in 1966; the non-trivial fact that it is unique up to isometry is due to Lusky \cite{lusky} in 1976.
An elementary proof has been recently found by Solecki and the second author~\cite{KS}. For a recent survey of the \Gurarii\ space and its non-separable versions we refer to \cite{kubis-gar}.

We shall construct a nonexpansive (i.e. of norm $\loe1$) linear operator $\map {\uop}\G \G$ with the following property: Given an arbitrary linear operator $\map T X Y$ between separable Banach spaces such that $\norm T \loe 1$, there exist isometric copies $X' \subs \G$ and $Y' \subs \G$ of $X$ and $Y$ respectively, such that $\img {\uop}{X'} \subs Y'$ and $\uop \rest X$ is isometric to $T$.
More formally, there exist isometric embeddings $\map i X \G$ and $\map j Y \G$ such that the following diagram is commutative.
$$\xymatrix{
\G \ar[r]^{\uop} & \G \\
X \ar[u]^i \ar[r]_T & Y \ar[u]_j
}$$
In other words, up to linear isometries, restrictions of $\uop$ to closed subspaces of $\G$ give \emph{all} nonexpansive linear operators between separable Banach spaces.

Furthermore, we show that the operator $\uop$ can be characterized by a condition similar to the one defining the \Gurarii\ space.

\section{Preliminaries}

We shall use standard notation concerning Banach space theory.
By $\nat$ we mean the set of all nonnegative integers.
We shall deal exclusively with nonexpansive linear operators, i.e., operators of norm $\loe 1$.
According to this agreement, a linear operator $\map fXY$ is an \emph{$\eps$-isometric embedding} if
$$ (1+\eps)^{-1} \cdot \norm {x} \leq  \norm {f(x)} \leq  \norm {x}$$
holds for every $x\in X$.
We shall often say ``\emph{$\eps$-embedding}'' instead of ``$\eps$-isometric embedding". In particular, an \emph{embedding} of one Banach space into another is a linear isometric embedding.
When dealing with a linear operator we shall always have in mind, besides its domain, also its \emph{co-domain}, which is just a (fixed in advance) Banach space containing the range (set of values) of the operator.

The \Gurarii\ space will be denoted by $\G$.

We shall use some standard category-theoretic notions.
Our basis is $\Ban$, the category of Banach spaces with linear operators of norm $\loe1$.
An important property of $\Ban$ is the following standard and well-known fact (see e.g. \cite{ACCGM}, \cite{gurarii} or \cite{pelczynski}).

\begin{lm}\label{faktone}
Let $\map i Z X$, $\map f Z Y$ be nonexpansive operators between Banach spaces.
Then there are nonexpansive operators $\map {g} X W$ and $\map {j} Y W$ such that
$$\xymatrix{
Y \ar[r]^{j} & W \\
Z \ar[r]_i \ar[u]^f & X \ar[u]_{g}
}$$
is a pushout square in $\Ban$.
Furthermore, if $i$ is an isometric embedding then so is $j$.
\end{lm}

It is worth mentioning the description of the pushout.
Namely, given $i$, $f$ as above, one usually defines $W = (X \oplus Y) \by \Delta$, 
where $X \oplus Y$ denotes the $\ell_1$-sum of $X$ and $Y$ and
$$\Delta = \setof{(i(z), -f(z))}{z \in Z}.$$
The operators $j$, $g$ are defined in the obvious way.

In case both $i$, $f$ are isometric embeddings, it can be easily seen that the unit ball of $W$ is the convex hull of the union of the unit balls of $X$ and $Y$, canonically embedded into $W$.
This remark will be used later.

\subsection{Correcting almost isometries}

Assume $\map f X Y$ is an $\eps$-embedding of Banach spaces, where $\eps > 0$.
It is natural to ask whether there exists an embedding $\map h X Y$ $\eps$-close to $f$.
Obviously, this may be impossible, since $Y$ may not contain isometric copies of $X$ at all.
Thus, a well-posed question is whether $f$ is $\eps$-close to some isometric embedding into some bigger Banach space containing $Y$.
This is indeed true, proved as Lemma~2.1 in \cite{KS}.
In fact, this is an elementary fact and very likely it appeared somewhere in the literature although the authors were unable to find it.
The proof of \cite[Lemma~2.1]{KS} uses linear functionals.
Below we provide a more elementary argument (coming from \cite{CsGwK}), at the same time showing that the ``correcting" isometric embedding is universal in the appropriate category.

Throughout this section we fix $\eps > 0$ and an $\eps$-embedding $\map f X Y$.
Actually, it is enough to require that $f$ satisfies
$$(1 - \eps) \norm x \loe \norm {f(x)} \loe (1 + \eps) \norm x,$$
although we consider nonexpansive operators only, therefore always $\norm{f(x)} \loe \norm x$.
Note that $1 - \eps < (1 + \eps)^{-1}$.
We define the following category $\fK(f,\eps)$.
The objects of $\fK(f,\eps)$ are pairs $\pair i j$ such that $\map i X Z$, $\map j Y Z$ are linear operators of norm $\loe 1$ such that
$$\norm{i - j \cmp f} \loe \eps.$$
Given two objects $a_1 = \pair {i_1}{j_1}$, $a_2 = \pair {i_2}{j_2}$, an arrow from $a_1$ to $a_2$ is a linear operator $h$ of norm $\loe 1$ such that
$$h \cmp i_1 = i_2 \oraz h \cmp j_1 = j_2.$$
By \cite[Lemma~2.1]{KS}, we know that if $f$ is an $\eps$-embedding then the category $\fK(f,\eps)$ contains an object $\pair i j$ such that both $i$ and $j$ are isometries.
Below we improve this fact, at least from the category-theoretic perspective.

\begin{lm}[cf. \cite{CsGwK}]\label{Lkluczowy}
The category $\fK(f,\eps)$ has an initial object $\pair {i_X}{j_Y}$ such that $\map {i_X} X {Z_0}$, $\map {j_Y} Y {Z_0}$ are isometries.

More precisely: $i_X$, $j_Y$ are canonical embeddings into $X \oplus Y$ endowed with the norm defined by the formula
$$\norm{v}_C = \inf\Bigsetof{\norm {x}_X + \norm {y}_Y + \eps \norm {w}_X}{v = \pair {x + w}{y - f(w)},\; x,w\in X,\; y\in Y},$$
where $\anorm_X$, $\anorm_Y$ are the norms of $X$ and $Y$ respectively.
\end{lm}

\begin{pf}
It is easy to check that $\anorm_C$ is indeed a norm.
In fact, the unit ball of $\anorm_C$ is the convex hull of the set $(\ubal X \times \sn 0) \cup (\sn 0 \times \ubal Y) \cup G$, where
$$G = \setof{\pair w{f(w)}}{\norm{w}_X \loe \eps^{-1}}.$$
Note that $\pair{i_X}{j_Y}$ is an object of $\fK(f,\eps)$, because $\norm{\pair w{-f(w)}} \loe \eps \norm {x}_X$ and $\norm{\pair x0}_C \loe \norm {x}_X$, $\norm{\pair 0y} \loe \norm {y}_Y$.

Fix an object $\pair ij$ of $\fK(f,\eps)$, and let $Z$ be the common range of $i$ and $j$.
Clearly, there exists a unique linear operator $\map h{X \oplus Y} Z$ such that $h \cmp i_X = i$ and $h \cmp j_Y = j$.
Namely, $h(x,y) = i(x) + j(y)$.
Note that $\norm{h(x,0)} \loe \norm{x}_X$, $\norm{h(0,y)} \loe \norm{y}_Y$, and $\norm{h(w,-f(w))} \loe \eps \norm{w}_X$.
The last inequality comes from the fact that $\pair i j$ is an object of $\fK(f,\eps)$.
It follows that $\norm{h(a)} \loe 1$ whenever $a$ is in the convex hull of $(\ubal X \times \sn 0) \cup (\sn 0 \times \ubal Y) \cup G$.
This shows that $\norm h \loe 1$, concluding the fact that $\pair {i_X}{j_Y}$ is an initial object of $\fK(f,\eps)$.

It remains to show that $i_X$ and $j_Y$ are isometries.
Fix $x\in X$.
Clearly, $\norm{\pair x0}_C \loe \norm {x}_X$.
On the other hand, for every $v\in Y$, $u,w\in X$ such that $u+w = x$ and $v - f(w) = 0$, we have
\begin{align*}
\norm {u}_X + \norm {v}_Y + \eps \norm {w}_X &= \norm {u}_X +\norm {f(w)}_Y + \eps \norm {w}_X\\
&\goe \norm {u}_X + (1 - \eps) \norm {w}_X + \eps \norm {w}_X\\
&\goe \norm {u+w}_X = \norm {x}_X.
\end{align*}
Passing to the infimum, we see that $\norm{\pair x0}_C \goe \norm {x}_X$.
This shows that $i_X$ is an isometric embedding.

Now fix $y\in Y$. Again, $\norm{\pair 0y}_C \loe \norm {y}_Y$ is clear.
Given $u,w\in X$, $v\in Y$ such that $u + w = 0$ and $v - f(w) = y$, we have
\begin{align*}
\norm {u}_X + \norm {v}_Y + \eps \norm {w}_X &= \norm {w}_X + \norm {v}_Y + \eps \norm {w}_X \\
&\goe \norm {v}_Y + \norm {f(w)}_Y\\
&\goe \norm {v - f(w)}_Y = \norm {y}_Y.
\end{align*}
Again, passing to the infimum we get $\norm{\pair 0y}_C \goe \norm {y}_Y$.
This shows that $j_Y$ is an isometric embedding and completes the proof.
\end{pf}

Note that when $\eps \goe 1$ then $f$ does not have to be an almost isometric embedding. In fact $f = 0$ can be taken into account. In such a case the initial object is just the coproduct $X \oplus Y$ with the $\ell_1$-norm.
In general, we shall denote by $X \oplus_{(f,\eps)} Y$ the space $X \oplus Y$ endowed with the norm described in Lemma~\ref{Lkluczowy} above.
Note that if $\almap f X Y$ is an $\eps$-embedding and $0 < \eps < \delta$ then $f$ is also a $\delta$-embedding, however the norm of $X \oplus_{(f,\eps)} Y$ is different from that of $X \oplus_{(f,\delta)} Y$.

The following statement will be used several times later.

\begin{lm}\label{Lkeycz}
Let $\eps, \delta > 0$ and let
$$\xymatrix{
X_0 \ar[d]_{T_0} \ar@{~>}[r]^{f_0} & Y_0 \ar[d]^{T_1} \\
X_1 \ar@{~>}[r]_{f_1} & Y_1
}$$
be a $\delta$-commutative diagram in $\Ban$ (i.e., $\norm{f_1 \cmp T_0 - T_1 \cmp f_0} \loe \delta$), such that $f_0$, $f_1$ are $\eps$-embeddings.
Then the operator
$$\map{ T_0 \oplus T_1 }{ X_0 \oplus_{(f_0,\eps+\delta)} Y_0 }{ X_1 \oplus_{(f_1,\eps)} Y_1 }$$
has norm $\loe 1$ and
\begin{equation}
(T_0 \oplus T_1) \cmp i_{X_0} = i_{X_1} \cmp T_0 \oraz
(T_0 \oplus T_1) \cmp j_{Y_0} = j_{Y_1} \cmp T_1.
\tag{$*$}\label{Eqtgno}
\end{equation}
\end{lm}

The situation is described in the following diagram, where the side squares are commutative, the bottom one is $\delta$-commutative, the left-hand side triangle $(\eps+\delta)$-commutative and the right-hand side triangle is $\eps$-commutative.
$$\xymatrix{
& X_0\oplus Y_0\ar@{->}[rrr]^{T_0\oplus T_1}& \ar@{-}[r]& \ar[r] &X_1\oplus Y_1\\
&\\
&Y_0 \ar@{^{(}->}[uu]|{j_{Y_0}} \ar@{->}[rrr]^{T_1}& \ar@{-}[r]& \ar[r] &Y_1 \ar@{^{(}->}[uu]|{j_{Y_1}} \\
X_0  \ar@{^{(}->}[uuur]|{i_{X_0}} \ar@{~>}[ur]^{f_0}  \ar@{->}[rrr]^{T_0}& \ar@{-}[r]& \ar[r] & X_1\ar@{~>}[ur]^{f_1} \ar@{^{(}->}[uuur]|{i_{X_1}}
}$$

\begin{pf}
By Lemma~\ref{Lkluczowy} applied to the category $\fK(f_0,\eps+\delta)$, there is a unique nonexpansive operator $\map S { X_0 \oplus_{(f_0,\eps+\delta)} Y_0 }{ X_1 \oplus_{(f_1,\eps)} Y_1 }$ satisfying (\ref{Eqtgno}) in place of $T_0 \oplus T_1$.
On the other, obviously $T_0 \oplus T_1$ satisfies (\ref{Eqtgno}), therefore $S = T_0 \oplus T_1$, showing that $T_0 \oplus T_1$ is nonexpansive. 
\end{pf}

\subsection{Rational operators}

We say that a Banach space $(X, \anorm)$ is \emph{rational} if $X$ is finite-dimensional and there exists a linear isomorphism $\map h {\Err^n}X$ such that
$$\norm x = \max_{i \loe k} |f_i(x)|, \qquad x \in X,$$
where $f_0,\dots,f_{k-1}$ are linear functionals preserving vectors with rational coordinates, that is, $\img {f_i h}{\Qyu^n} = \Qyu$ for $i<k$.
Very formally, a rational Banach space is a triple of the form $(X, \anorm, h)$, where $(X, \anorm)$ and $h$ are as above.
This notion is needed for catching countably many spaces that approximate the class of all finite-dimensional Banach spaces.
In other words, a finite-dimensional Banach space is rational if there is a coordinate-wise system (induced by a linear isomorphism from $\Err^n$ and by the standard basis of $\Err^n$) such that the closed unit ball is the convex hull of finitely many vectors, each of them having rational coordinates.
Note that ``being rational" depends both on the norm and on the coordinate-wise system.
For instance, the two-dimensional Hilbert space is not rational, however the space $\Err^2$ endowed with a scaled $\ell_1$-norm 
$$\norm{(x,y)} = \sqrt 2 (|x| + |y|)$$
is rational, which is witnessed by the isomorphism $h(v) = \sqrt 2 v$, $v \in \Err^2$.
Later on, forgetting this example, when considering $\Err^n$ as a rational Banach space, we shall always have in mind the usual coordinate-wise system.

We shall also need the notion of a rational operator.
Namely, an operator $\map T {X}{Y}$ is rational if $\norm T \loe 1$ and $\img {T h}{\Qyu^m} \subs \img g {\Qyu^n}$, where 
$\map h {\Err^m} X$ and $\map g {\Err^n} Y$ are linear isomorphisms with respect to which $X$ and $Y$ are rational Banach spaces.
Note that there are, up to isometry, only countably many rational operators.
Note also that being a rational operator again depends on fixed linear isomorphism inducing coordinate-wise systems.

\begin{lm}\label{LmNormyRacjonalne}
Let $X \subs Y$ be finite-dimensional Banach spaces and assume that $Y = \Err^m$ so that $X$ is its rational subspace and the norm $\anorm_Y$ is rational when restricted to $X$.
Then for every $\delta>0$ there exists a rational norm $\anorm_Y'$ on $Y$ that is $\delta$-equivalent to $\anorm_Y$ and such that $\norm{x}_Y = \norm{x}'_Y$ for every $x\in X$.
\end{lm}

\begin{pf}
Let $\Phi$ be a finite collection of rational functionals on $X$ such that
$$\norm{x}_X = \max_{\phi \in \Phi}|\phi(x)|$$ for every $x \in X$.
By the Hahn-Banach Theorem, we may assume that each $\phi \in \Phi$ is actually a rational functional on $Y$ that has ``almost'' the same norm as its restriction to $X$.
Finally, enlarge $\Phi$ to a finite collection $\Phi'$ by adding finitely many rational functionals so that the new norm induced by $\Phi'$ will become $\delta$-equivalent to $\anorm_Y$.
\end{pf}

\begin{lm}\label{Lmkorektwym}
Let $\map {T_0} {X_0} {Y_0}$ be a rational operator, $\eps>0$, and let $\map T X Y$ be an operator of norm $\loe1$ extending $T_0$ and such that $X \sups X_0$, $Y \sups Y_0$ are finite-dimensional. Let $\anorm_X$, $\anorm_Y$ denote the norms of $X$, $Y$.
Then there exist rational norms $\anorm_X'$ and $\anorm_Y'$ on $X$ and $Y$, respectively, such that
\begin{enumerate}
	\item[(i)] $T$ is a rational operator from $(X,\anorm_X')$ to $(Y,\anorm_Y')$,
	\item[(ii)] $\anorm_X'$ is $\eps$-equivalent to $\anorm_X$ and $\anorm_Y'$ is $\eps$-equivalent to $\anorm_Y$,
	\item[(iii)] $X_0 \subs X$ and $Y_0 \subs Y$ are rational isometric embeddings when $X$, $Y$ are endowed with $\anorm_X'$, $\anorm_Y'$ and $X_0$, $Y_0$ are endowed with their original norms.
\end{enumerate}
\end{lm}

\begin{pf}
For simplicity, let us assume that $X = X_0 \oplus \Err u$ and $Y = Y_0 \oplus \Err v$ and either $T(u)=v$ or $T(u)$ is a rational vector in $X_0$.
The general case will follow by induction.
Let $\map {h_0}{\Err^m}{X_0}$, $\map {g_0}{\Err^n}{Y_0}$ be linear isomorphisms witnessing that $X_0$, $Y_0$ are rational Banach spaces and that $T_0$ is a rational operator.
Extend $h_0$, $g_0$ to $\map h{\Err^{m+1}}X$, $\map g{\Err^{n+1}}Y$ by setting $h(e_{m+1}) = u$, $g(e_{n+1}) = v$, where $e_i$ denotes the $i$th vector from the standard vector basis of $\Err^k$ ($k \goe i$).
Note that $T$ will become a rational operator, as long as we define suitable rational norms on $X$ and $Y$.

Fix $\delta>0$.
Let $\anorm_X'$ and $\anorm_Y'$ be obtained from Lemma~\ref{LmNormyRacjonalne}.
Conditions (i) and (iii) are obviously satisfied.
The only obstacle is that the operator $T$ may not be nonexpansive with respect to these new norms.
However, we have $\norm{T x}'_Y \loe (1+\delta) \norm {T x}_Y \loe (1+\delta)\norm {x}_X \loe (1+\delta)^2 \norm {x}'_X$.
Assuming that $\delta$ is rational, we can replace $\anorm'_X$ by $(1+\delta)^2 \anorm'_X$, so that $T$ is again nonexpansive.
Finally, if $\delta$ is small enough, then condition (ii) holds.
\end{pf}

\subsection{The \Gurarii\ property}

Before we construct the isometrically universal operator, we consider its crucial property which is similar to the condition defining the \Gurarii\ space.
Namely, we shall say that a linear operator $\map \Omega U V$ has the \emph{\Gurarii\ property} if $\norm \Omega \loe1$ and the following condition is satisfied.
\begin{enumerate}
	\item[($G$)] Given $\eps>0$, given a nonexpansive operator $\map {T}{X}{Y}$ between finite-dimensional spaces, given $X_0 \subs X$, $Y_0 \subs Y$ and isometric embeddings $\map i {X_0}U$, $\map j {Y_0}V$ such that $\Omega \cmp i = j \cmp (T\rest X_0)$, there exist $\eps$-embeddings $\map {i'} X U$, $\map {j'} Y V$ satisfying
$$\norm{i'\rest X_0 - i}<\eps, \quad \norm{j'\rest Y_0 - j}<\eps, \oraz \norm{\Omega \cmp i'  - j' \cmp T}<\eps.$$
\end{enumerate}
We shall also consider condition ($G^*$) which is, by definition, the same as ($G$) with the stronger requirement that $\Omega \cmp i' = j' \cmp T$.
We shall see later that ($G$) is equivalent to ($G^*$).
In the next section we show that an operator with the \Gurarii\ property exists.

\section{The construction}

Fix two real vector spaces $U$, $V$, each having a fixed countable infinite vector basis, which provides a coordinate system and the notion of rational vectors (namely, rational combinations of the vectors from the basis).
For the sake of brevity, we may assume that $U$,$V$ are disjoint, although this is not essential.
We shall construct rational subspaces of $U$ and $V$.
Notice the following trivial fact: given a rational space $X \subs U$ (that is, $X$ is finite-dimensional and its norm is rational in $U$), given a rational isometric embedding $\map e X Y$ (so $Y$ is also a rational space), there exists a rational extension $X'$ of $X$ in $U$ and a bijective rational isometry $\map h Y {X'}$ such that $h \cmp e$ is identity on $X$.
In other words, informally, every rational extension of a rational space ``living" in $U$ is realized in $U$.
Of course, the same applies to $V$.

We shall now construct a sequence of rational operators $\map {F_n}{U_n}{V_n}$ such that
\begin{enumerate}
	\item[(a)] $U_n$ is a rational subspace of $U$ and $V_n$ is a rational subspace of $V$.
	\item[(b)] $F_{n+1}$ extends $F_n$ (in particular, $U_n \subs U_{n+1}$ and $V_n \subs V_{n+1}$).
	\item[(c)] Given $\Ntr$, given rational embeddings $\map i{U_n}X$, $\map j{V_n}Y$ and given a rational operator $\map T X Y$ such that $T \cmp i = j \cmp F_n$, there exist $m > n$ and rational embeddings $\map {i'}X{U_m}$, $\map {j'}Y{V_m}$ satisfying $j' \cmp T = F_m \cmp i'$ and such that $i' \cmp i$ and $j' \cmp j$ are identities on $U_n$ and $V_n$, respectively.
\end{enumerate}
For this aim, let $\Ef$ denote the family of all triples $\triple T e k$, where $\map T X Y$ is a rational operator, $k$ is a natural number and $e = \pair i j$ is a pair of rational embeddings like in condition (c) above, namely $\map {i}{X_0}X$, $\map {j}{Y_0}Y$, where $X_0$ and $Y_0$ are rational subspaces of $U$ and $V$, respectively.
We also assume that $X$, $Y$ are rational subspaces of $U$, $V$, therefore the family $\Ef$ is indeed countable.
Enumerate it as $\sett{\triple {T_n}{e_n}{k_n}}{\ntr}$ so that each $\triple T e k$ appears infinitely many times.

We now start with $U_0=0$, $V_0=0$ and $F_0=0$.
Fix $n > 0$ and suppose $\map {F_{n-1}}{U_{n-1}}{V_{n-1}}$ has been defined.
We look at triple $\triple {T_n}{e_n}{k_n}$ and consider the following condition, where $e_n = \pair {i_n}{j_n}$:
\begin{enumerate}
	\item[($*$)] $k_n = k < n$, the domain of $i_n$ is $U_k$, the domain of $j_n$ is $V_k$ and $T \cmp i_n = j_n \cmp F_k$.
\end{enumerate}
If ($*$) fails, we set $F_n = F_{n-1}$ (and $U_n = U_{n-1}$, $V_n = V_{n-1}$).
Suppose now that ($*$) holds and $\map{i_n}{U_k}X$, $\map{j_n}{V_k}Y$.
Using the push-out property (more precisely, its version for rational operators), we find rational spaces $X' \sups U_{n-1}$ and $Y' \sups V_{n-1}$ (and the inclusions are rational embeddings), together with rational embeddings $\map {i'}X{X'}$, $\map {j'}Y{Y'}$ such that $i' \cmp i_n$ is identity on $U_k$ and $j' \cmp j_n$ is identity on $V_k$.
As we have mentioned, we may ``realize'' the spaces $X'$ and $Y'$ inside $U$ and $V$, respectively.
We set $U_n = X'$, $V_n = Y'$ and we define $F_n$ to be the unique operator from $U_n$ to $V_n$ obtained from the push-outs.
In particular, $F_n$ extends $F_{n-1}$ and satisfies $j' \cmp T = F_n \cmp i'$.
This completes the description of the construction of $\sett{F_n}{\ntr}$.

Note that the sequence satisfies (a)--(c).
Only condition (c) requires an argument.
Namely, fix $n$ and $T,i,j$ as in (c). Find $m>n$ such that $\triple {T_m}{e_m}{k_m} = \triple T e n$, where $e = \pair i j$.
Then at the $m$th stage of the construction condition ($*$) is fulfilled and therefore $F_m$ witnesses that (c) holds.

Now denote by $U_\infty$ and $V_\infty$ the completions of $\bigcup_{\Ntr}U_n$ and $\bigcup_{\Ntr}V_n$ and
let $\map{F_\infty}{U_\infty}{V_\infty}$ the unique extension of $\bigcup_{\Ntr}F_n$.

\begin{prop}
The operator $F_\infty$ satisfies condition ($G^*$) and therefore has the \Gurarii\ property.
\end{prop}

\begin{pf}
Fix finite-dimensional Banach spaces $X_0 \subs X_1$, $Y_0 \subs Y_1$, fix isometric embeddings $\map i {X_0}{U_\infty}$, $\map j {Y_0}{V_\infty}$ so that $F_\infty \cmp i = j \cmp T_0$.
Furthermore, fix two non\-ex\-pan\-sive operators $\map {T_i}{X_i}{Y_i}$ for $i=0,1$ such that $T_1$ extends $T_0$.
Fix $\eps>0$.

We need to find $\eps$-isometric embeddings $\map f{X_1}{U_\infty}$, $\map g{Y_1}{V_\infty}$ such that 
$$\norm{f\rest X_0 - i}<\eps, \quad \norm{g\rest Y_0 - j}<\eps,$$
and $F_\infty \cmp f = g \cmp T_1$.
Let $\delta = \eps/3$.
The remaining part of the proof is divided into four steps:

\paragraph{Step 1}
We first ``distort" the embeddings $i$, $j$, so that their images will be some $U_n$ and $V_n$, respectively.
Formally, we find $\delta$-isometric embeddings $\map {i_0}{X_0}{U_n}$, $\map {j_0}{Y_0}{V_n}$ for some fixed $n$, such that $\norm {i-i_0} < \delta$, $\norm{j-j_0}<\delta$ and $$\norm{j_0 \cmp T_0 - F_n \cmp i_0} < \delta.$$

\paragraph{Step 2}
Applying Lemma~\ref{Lkeycz} for obtaining finite-dimensional spaces $X_2 \sups X_0$, $Y_2 \sups Y_0$ with isometric embeddings $\map k {U_n}{X_2}$, $\map \ell {V_n}{Y_2}$, together with a nonexpansive operator $\map {T_2}{X_2}{Y_2}$ extending $T_0$ and satisfying $T_2 \cmp k = \ell \cmp F_n$.

\paragraph{Step 3}
We now use the Pushout Lemma for two pairs of embeddings: $X_0 \subs X_1$, $X_0 \subs X_2$ and $Y_0 \subs Y_1$, $Y_0 \subs Y_2$; we obtain a further extension of the operator $T_2$.
Thus, in order to avoid too many objects, we shall assume that $X_1 \subs X_2$, $Y_1 \subs Y_2$ and the operator $T_2$ extends both $T_0$ and $T_1$.
At this point, we may actually forget about $T_1$, replacing it by $T_2$.

\paragraph{Step 4}
Apply Lemma~\ref{Lmkorektwym} in order to change the norms of $X_2$ and $Y_2$ by $\delta$-equivalent ones, so that $T_2$ becomes a rational operator extending $F_n$.
Using (c), we can now ``realize'' $T_2$ in $F_m$ for some $m > n$.
Formally, there are isometric embeddings $\map {i_2}{X_2}{U_m}$, $\map {j_2}{Y_2}{V_m}$ satisfying $F_m \cmp i_2 = j_2 \cmp T_2$.
Coming back to the original norms of $X_2$ and $Y_2$, we see that the embeddings $i_2$ and $j_2$ are $\delta$-isometric and 
their restrictions to $X_0$ and $Y_0$ are $(2\delta)$-close to $i_0$ and $j_0$, respectively, hence $\eps$-close to $i$ and $j$ (recall that $\delta = \eps/3$).
This completes the proof.
\end{pf}

\section{Properties}

We now show that the domain and the co-domain of an operator with the \Gurarii\ property acting between separable spaces is the \Gurarii\ space (thus, justifying the name) and later we show its uniqueness as well as some kind of homogeneity.

\subsection{Recognizing the domain and the co-domain}

Recall that a separable Banach space $W$ is linearly isometric to the \Gurarii\ space $\G$ if and only if it satisfies the following condition:
\begin{enumerate}
	\item[($\gurdef$)] Given finite-dimensional spaces $X_0 \subs X$, given $\eps>0$, given an isometric embedding $\map i {X_0}W$, there exists an $\eps$-embedding $\map f X W$ such that $\norm{f \rest X_0 - i} \loe \eps$.
\end{enumerate}
Usually, the condition defining the \Gurarii\ space is stronger, namely, it is required that $f \rest X_0 = i$.
For our purposes, the formally weaker condition ($\gurdef$) is more suitable.
It is not hard to see that both conditions are actually equivalent, see~\cite{kubis-gar} for more details.

\begin{tw}\label{Thmrigbiw}
Let $\map \Omega U V$ be a linear operator with the \Gurarii\ property, where $U,V$ are separable.
Then both $U$ and $V$ are linearly isometric to the \Gurarii\ space.
\end{tw}

\begin{pf}
(1) $U$ is isometric to $\G$.

\noindent
Fix finite-dimensional spaces $X_0 \subs X$ and fix an isometric embedding $\map i {X_0}U$. Fix $\eps>0$.
Let $Y_0 = \img \Omega {\img i {X_0}}$ and let $T_0 = \Omega \rest \img i{X_0}$, treated as an operator into $Y_0$.
Applying the pushout property (Lemma~\ref{faktone}), we find a finite-dimensional space $Y \sups Y_0$ and a nonexpansive linear operator $\map T X Y$ extending $T_0$.
Applying condition ($G$), we get in particular an $\eps$-embedding $\map f X U$ satisfying $\norm{f\rest X_0 - i} \loe \eps$.
This shows that $U$ satisfies ($\gurdef$).

(2) $V$ is isometric to $\G$.

\noindent
Fix $\eps>0$ and fix finite-dimensional spaces $Y_0 \subs Y$ and an isometric embedding $\map j {Y_0} V$.
Let $X_0 = \sn0 = X$ and let $\map {T_0}{X_0}{Y_0}$, $\map T X Y$ be the 0-operators.
Applying ($G$), we get an $\eps$-embedding $\map {j'} Y V$ satisfying $\norm{j' \rest Y_0 - j}\loe \eps$, showing that $V$ satisfies ($\gurdef$).
\end{pf}

From now on, we shall denote by  $\uop$ the operator constructed in the previous section. According to the results above, this operator has the \Gurarii\ property and it is of the form $\map \uop \G \G$.

\subsection{Universality}

We shall now simplify the notation, in order to avoid too many parameters and shorten some arguments.
Namely, given nonexpansive linear operators $\map S X Y$, $\map T Z W$,
a pair $i = \pair {i_0}{i_1}$ of isometric embeddings of the form $\map {i_0} X Z$, $\map {i_1} Y W$ and satisfying $T \cmp i_0 = i_1 \cmp S$, will be called an \emph{embedding of operators} from $S$ into $T$ and we shall write $\map i S T$.
Now fix $\eps > 0$ and suppose that $f = \pair {f_0}{f_1}$ is a pair of $\eps$-embeddings of the form $\map {f_0}X Z$, $\map {f_1}Y W$ satisfying $\norm{T \cmp f_0 - f_1 \cmp S} \loe \eps$.
We shall say that $f$ is \emph{$\eps$-embedding of operators} from $S$ into $T$ and we shall write $\almap f S T$.
Finally, an \emph{almost embedding of operators} of $S$ into $T$ will be, by definition, an $\eps$-embedding of $S$ into $T$ for some $\eps>0$.
The composition of (almost) embeddings of operators is defined in the obvious way.
Given two almost embeddings of operators $\map f S T$, $\map g S T$, we shall say that $f$ is \emph{$\eps$-close} to $g$ (or that $f$, $g$ are \emph{$\eps$-close}) if $\norm{f_0 - g_0} \loe \eps$ and $\norm{f_1 - g_1} \loe \eps$, where $f = \pair {f_0}{f_1}$ and $g = \pair {g_0}{g_1}$.

The notation described above is in accordance with category-theoretic philosophy: almost embeddings of operators obviously form a category, which is actually a special case of much more general constructions on categories, where the objects are diagrams of certain shape.

Now, observe that a consequence of the pushout property stated in Lemma~\ref{faktone} (that we have already used) says that for every two embeddings of operators $\map i S T$, $\map j S R$ there exist embeddings of operators $\map {i'} T P$, $\map {j'} R P$ such that $i' \cmp i = j' \cmp j$.
The crucial property of almost embeddings is Lemma~\ref{Lkeycz} which says, in the new notation, that for every $\eps$-embedding of operators $\map f S T$ there exist embeddings of operators $\map i S R$, $\map j T R$ such that $j \cmp f$ is $(2\eps)$-close to $i$.

Before proving the universality of $\uop$, we formulate a statement which is crucial for the proof.

\begin{lm}\label{Lmrgeihei}
Let $\Omega$ be an operator with the \Gurarii\ property.
Assume $\eps>0$ and $\almap f T \Omega$ is an $\eps$-embedding of operators and $\map j T R$ is an embedding of operators, where both $T$ and $R$ act between finite-dimensional spaces.
Then for every $\delta > 0$ there exists a $\delta$-embedding of operators $\almap g R \Omega$ whose composition with $j$ is $(2\eps+\delta)$-close to $f$.
\end{lm}

\begin{pf}
We first replace $f$ by $\almap {f_1} T {T_1}$, where $T_1$ is some restriction of $\Omega$ to finite-dimensional spaces (so $f = e \cmp f_1$, where the components of $e$ are inclusions).

Using Lemma~\ref{Lkeycz}, we find isometries of operators $\map i {T_1} S$, $\map {j_1} T S$, where $S$ is an operator between finite-dimensional spaces and $j_1$ is $(2\eps)$-close to $i \cmp f_1$.
Applying the amalgamation property to $j_1$ and $j$, we find embeddings of operators $\map k R {\tilde S}$ and $\map {j_2} S {\tilde S}$ satisfying $j_2 \cmp j_1 = k \cmp j$.
In order to avoid too many parameters, we replace $S$ by $\tilde S$ and $j_1$ by $j_2 \cmp j_1$.
By this way, we have an embedding of operators $\map k R S$ such that $k \cmp j = j_1$ and still $j_1$ is $(2\eps)$-close to $i \cmp f_1$.

Now, condition ($G$) in our terminology says that there is a $\delta$-embedding of operators $\almap \ell S \Omega$ whose composition with $i$ is $\delta$-close to the inclusion $\map e {T_1}\Omega$.
Finally, $g = \ell \cmp k$ is the required $\delta$-embedding, because $g \cmp j$ is $(2\eps + \delta)$-close to $f$.
\end{pf}

\begin{tw}
Given a nonexpansive linear operator $\map T X Y$ between separable Banach spaces, there exist isometric embeddings $\map i X \G$, $\map j Y \G$ such that $\uop \cmp i = j \cmp \uop$, that is, the following diagram is commutative.
$$\xymatrix{
\G \ar[r]^{\uop} & \G \\
X \ar[u]^i \ar[r]_T & Y \ar[u]_j
}$$
\end{tw}

\begin{pf}
We first ``decompose'' $T$ into a chain $T_0 \subs T_1 \subs T_2 \subs \cdots$ so that $\map {T_n}{X_n}{Y_n}$ and $X_n$, $Y_n$ are finite-dimensional spaces.
Formally, we construct inductively two chains of finite-dimensional spaces $\sett{X_n}{\ntr}$, $\sett{Y_n}{\ntr}$ such that $\img T {X_n} \subs Y_n$ and $\bigcup_{\ntr}X_n$ is dense in $X$, and $\bigcup_{\ntr}Y_n$ is dense in $Y$.
It is clear that such a decomposition is always possible and the operator $T$ is determined by the chain $\sett{T_n}{\ntr}$.

Let $\eps_n = 2^{-n}$.
We shall construct almost embeddings of operators $\map {i_n} {T_n} \uop$ so that the following conditions are satisfied:
\begin{enumerate}
	\item[(i)] $i_n$ is an $\eps_n$-embedding of $T_n$ into $\uop$.
	\item[(ii)] $i_{n+1} \rest T_n$ is $(3\eps_n)$-close to $i_n$.
\end{enumerate}
Once we assume that $i_0$ is the 0-operator between the 0-subspaces of $X$ and $Y$, there is no problem to start the inductive construction.
Fix $n \goe 0$ and suppose $i_n$ has already been constructed.
Let $Z = \img {T_n}{X_n}$.
By Lemma~\ref{Lmrgeihei} applied to $i_n$ with $\eps=\delta=\eps_n$, we find $i_{n+1}$ satisfying (ii).
Thus, the construction can be carried out.

Finally, there is a unique operator $\map {i_\infty} T \G$ that extends all the $i_n$s.
Formally, both components of $i_\infty$ are uniquely determined by the completion of the pointwise limit of the sequence $\sett {i_n}{\ntr}$.
The two components of $i_\infty$ are the required isometric embeddings.
\end{pf}

\subsection{Uniqueness and almost homogeneity}

We start with the main, somewhat technical, lemma from which we easily derive all the announced properties of $\uop$.
We shall say that an operator $\map f X Y$ is a \emph{strict} $\eps$-embedding if it is a $\delta$-isometric embedding for some $0<\delta<\eps$.

\begin{lm}\label{LmMejnn}
Assume $\map \Omega U V$ and $\map {\Omega'}{U'}{V'}$ are two operators between separable Banach spaces, both with the \Gurarii\ property (that is, nonexpansive and with property ($G$)).
Assume $X_0 \subs U$, $Y_0 \subs V$, $X'_0 \subs U'$ and $Y'_0 \subs V'$ are finite-dimensional spaces.
Fix $\eps > 0$ and let $T_0 = \Omega \rest X_0$, $T_0' = \Omega' \rest X'_0$.
Assume further that $\map {i_0}{X_0}{X'_0}$ and $\map {j_0}{Y_0}{Y'_0}$ are strict $\eps$-embeddings satisfying
$$\norm{T_0' \cmp i_0 - j_0 \cmp T_0} < \eps.$$
Then there exist bijective linear isometries $\map I U{U'}$ and $\map J V{V'}$ such that
$$J \cmp \Omega = \Omega' \cmp I$$
and $\norm{I\rest X_0 - i_0} < \eps$, $\norm{J\rest Y_0 - j_0} < \eps$.
\end{lm}

Before proving this lemma, we formulate and prove some of its consequences.
First of all, let us say that an operator $\map \Omega U V$ is \emph{almost homogeneous} if the following condition is satisfied.
\begin{enumerate}
	\item[(AH)] Given $\eps>0$, given finite-dimensional spaces $X_0, X_1 \subs U$, $Y_0, Y_1 \subs V$ such that $\img \Omega {X_0} \subs Y_0$, $\img \Omega {X_1} \subs Y_1$, given linear isometries $\map i{X_0}{X_1}$, $\map j{Y_0}{Y_1}$ such that $\Omega \cmp i = j \cmp \Omega$, there exist bijective linear isometries $\map I U U$, $\map J V V$ such that $$\Omega \cmp I = J \cmp \Omega$$ and $\norm{I \rest X_0 - i} \loe \eps$, $\norm{J \rest Y_0 - j} \loe \eps$.
\end{enumerate}
Eliminating $\eps$ from this definition, we obtain the notion of a \emph{homogeneous operator}. We shall see in a moment, using our knowledge on the \Gurarii\ space, that no operator between separable Banach spaces can be homogeneous.

\begin{tw}\label{Thmenrgog}
Let $\Omega$ be a linear operator with the \Gurarii\ property, acting between separable Banach spaces.
Then $\Omega$ is isometric to $\uop$ in the sense that there exist bijective linear isometries $I$, $J$ such that $\uop \cmp I = J \cmp \Omega$.

Furthermore, $\Omega$ is almost homogeneous.
\end{tw}

\begin{pf}
We already know that $\uop$ has the \Gurarii\ property.
Applying Lemma~\ref{LmMejnn} to the zero operator, we obtain the required isometries $I$, $J$.

In order to show almost homogeneity, apply Lemma~\ref{LmMejnn} again to the operator $\Omega$ (on both sides) and to the embeddings $i$, $j$ specified in condition (AH).
\end{pf}

Let us say that an operator $\Omega$ is \emph{isometrically universal} if, up to isometries, its restrictions to closed subspaces provide all operators between separable Banach spaces whose norms do not exceed the norm of $\Omega$.

\begin{uwgi}
No bounded linear operator between separable Banach spaces can be isometrically universal and homogeneous.
\end{uwgi}

\begin{pf}
Suppose $\Omega$ is such an operator and consider $G = \ker \Omega$.
Then $G$ would have the following property: Every isometry between finite-dimensional subspaces of $G$ extends to an isometry of $G$.
Furthermore, $G$ contains isometric copies of all separable spaces, because $\Omega$ is assumed to be isometrically universal.
On the other hand, it is well-known that no separable Banach space can be homogeneous and isometrically universal for all finite-dimensional spaces, since this would be necessarily the \Gurarii\ space, which is not homogeneous (see \cite{gurarii} or \cite{kubis-gar}).
\end{pf}

It remains to prove Lemma~\ref{LmMejnn}.
It will be based on the ``approximate back-and-forth argument", similar to the one in \cite{KS}.
In the inductive step we shall use the following fact, formulated in terms of almost embeddings of operators.

\begin{claim}\label{Clejpwrt}
Assume $\Omega$ is an operator with the \Gurarii\ property, $\eps > 0$ and $\almap f T R$ is an $\eps$-embedding of operators acting on finite-dimensional spaces, and $\map e T \Omega$ is an embedding of operators.
Then for every $\delta > 0$ there exists a $\delta$-embedding $\almap g R \Omega$ such that $g \cmp f$ is $(2\eps+\delta)$-close to $e$.
\end{claim}

\begin{pf}
Using Lemma~\ref{Lkeycz}, we find embeddings of operators $\map i T S$, $\map j R S$ such that $j \cmp f$ is $(2\eps)$-close to $i$.
Property ($G$) tells us that there exists a $\delta$-embedding $\almap h S \Omega$ such that $h \cmp i$ is $\delta$-close to $e$.
Finally, $g = h \cmp j$ is $(2\eps+\delta)$-close to $e$.
\end{pf}

The usefulness of the above claim comes from the fact that $\delta$ can be arbitrarily small comparing to $\eps$.

\begin{pf}[Proof of Lemma~\ref{LmMejnn}]
We first choose $0 < \eps_0 < \eps$ such that $k_0 := \pair{i_0}{j_0}$ is an $\eps_0$-embedding of operators.

Our aim is to build two sequences $\ciag k$ and $\ciag \ell$ of almost embeddings of operators between finite subspaces of $U$, $V$, $U'$, $V'$.
Notice that once we are given operators and almost embeddings as in the statement of Lemma~\ref{LmMejnn}, we are always allowed to enlarge the co-domains (namely the spaces $Y_0$ and $Y'_0$) to arbitrarily big finite-dimensional subspaces of $V$ and $V'$, respectively.
This is important for showing that our sequences of almost embeddings will ``converge'' to bijective isometries.

The formal requirements are as follows.
We choose sequences $\ciag u$, $\ciag v$, $\sett{u'_n}{\ntr}$, $\sett{v'_n}{\ntr}$ that are linearly dense in $U$, $V$, $U'$, $V'$, respectively.
We fix a decreasing sequence $\ciag \eps$ of positive real numbers satisfying
\begin{equation}
3\sum_{n=1}^\infty \eps_n < \eps - \eps_0,
\tag{s}\label{EqSerious}
\end{equation}
where $\eps_0<\eps$ is as above.
We require that:
\begin{enumerate}
	\item[(1)] $\almap {k_n}{T_n}{T'_n}$, $\almap {\ell_n}{T'_n}{T_{n+1}}$ are almost embeddings of operators; $T_n$ is a restriction of $\Omega$ to some pair of finite-dimensional spaces and $T'_n$ is a restriction of $\Omega'$ to some pair of finite-dimensional spaces.
	\item[(2)] $k_n$ is an $\eps_n$-embedding of operators and $\ell_n$ is an $\eps_{n+1}$-embedding of operators.
	\item[(3)] The composition $\ell_n \cmp k_n$ is $(2\eps_n+\eps_{n+1})$-close to the identity (formally, to the inclusion $T_n \subs T_{n+1}$).
	\item[(4)] $k_{n} \cmp \ell_{n-1}$ is $(2\eps_{n}+\eps_{n+1})$-close to the identity.
	\item[(5)] $u_n$ belongs to the domain of $T_{n+1}$, $v_n$ belongs to its co-domain; similarly for $u_n'$, $v_n'$ and $T'_{n+1}$.
\end{enumerate}
Fix $n \goe 0$ and suppose that $k_n$ and $\ell_{n-1}$ have already been constructed (if $n=0$ then we ignore condition (4)).
We apply Claim~\ref{Clejpwrt} twice: first time to $k_n$ with $\eps = \eps_n$ and $\delta = \eps_{n+1}$, thus obtaining $\ell_n$; second time to $\ell_n$ with $\eps = \eps_{n+1}$ and $\delta = \eps_{n+2}$, thus obtaining $k_{n+1}$.
Between these two steps, we choose a sufficiently big operator $T_{n+1}$ which is a restriction of $\Omega$ to some finite-dimensional spaces.
Also, after obtaining $k_{n+1}$, we choose a sufficiently big operator $T'_{n+1}$ contained in $\Omega'$ and acting between finite-dimensional spaces.
By this way, we may ensure that condition (5) holds.

Fix $n>0$. We have the following (non-commutative) diagram of almost embeddings of operators
$$\xymatrix{
T_n \ar@{~>}[rrr]^{k_n} & & & T'_n \\
T_{n-1} \ar[u] \ar@{~>}[rrr]_{k_{n-1}} & & & T'_{n-1} \ar[u] \ar@{~>}[lllu]_{\ell_{n-1}}
}$$
in which the vertical arrows are inclusions, the lower triangle is $(2\eps_{n-1}+\eps_n)$-com\-mu\-ta\-tive by (3), and the upper triangle is $(2\eps_n+\eps_{n+1})$-commutative by (4).
Formally, these relations are true for the two components of all ``arrows" appearing in this diagram.

Using the triangle inequality of the norm and the fact that all operators are nonexpansive, we conclude that $k_n$ restricted to $T_{n-1}$ is $\eta_n$-close to $k_{n-1}$, where
$\eta_n = 2\eps_{n-1} + 3\eps_n + \eps_{n+1}$.
In particular, both components of the sequence $\ciag k$ are point-wise convergent and by (2) the completion of the limit defines an isometric embedding of operators $\map K \Omega {\Omega'}$.

Interchanging the roles of $\Omega$ and $\Omega'$, we deduce that the sequence $\ciag \ell$ converges to an isometric embedding of operators $\map L {\Omega'}\Omega$. Conditions (3),(4) say that $L$ is the inverse of $K$, therefore $K$ is bijective.
Finally, denoting $K = \pair I J$, we see that $I$, $J$ are as required, because $\sum_{n=1}^\infty \eta_n < 2\eps_0 + 6\sum_{n=1}^\infty\eps_n < 2\eps$.
\end{pf}

\subsection{Kernel and range}

We finally show some structural properties of our operator.

\begin{tw} The operator $\uop$ is surjective and its kernel is linearly isometric to the \Gurarii\ space.
\end{tw}

\begin{pf}
We first show that $\ker \uop$ is isometric to $\G$.

Fix finite-dimensional spaces $X_0 \subs X$ and fix an isometric embedding $\map i{X_0}{\ker \Omega}$ and fix $\eps > 0$.
Let $Y_0 = \sn0 = Y$ and let $\map j{Y_0}Y$ be the 0-operator.
Applying ($G^*$), we obtain $\eps$-embeddings $\map {i'} X U$ and $\map {i'} Y V$ such that $\norm{i'\rest X_0 - i} \loe \eps$, $\norm{j'\rest Y_0 - j} \loe \eps$ and $j' \cmp 0 = \Omega \cmp i'$, where $0$ denotes the 0-operator from $X$ to $Y$.
By the last equality, $i'$ maps $X$ into $\ker \Omega$, showing that $\ker \Omega$ satisfies ($\gurdef$).

In order to show that $\img \uop \G = \G$, fix $v \in \G$ and consider the zero operator $\map {T_0}{X_0}{Y_0}$, where $X_0 = \sn 0$ and $Y_0$ is the one-dimensional subspace of $\G$ spanned by $v$.
Let $X = Y = Y_0$ and let $\map T X Y$.
Let $\map i {X_0}\G$ and $\map j {Y_0}\G$ denote the inclusions.
Applying condition ($G^*$) with $\eps=1$, we obtain linear operators $\map {i'} X \G$, $\map {j'} Y \G$ such that $\norm{i'\rest X_0 - i}<1$, $\norm{j'\rest Y_0 - j}<1$ and $\uop \cmp i' = j' \cmp T$.
Obviously, $j'=j$ and hence $\uop(i'(v))= T(v) = v$.
\end{pf}

\subsection*{Acknowledgments}
The authors would like to thank Vladimir M\"uller for pointing out references concerning universal operators on Hilbert spaces.
The second author would like to thank the Hausdorff Research Institute for Mathematics (Bonn, September 2013) for their support and warm hospitality.


\begin{thebibliography}{9}

\bibitem{AmbMul} C. Ambrozie, V. M\"uller,
{\it Commutative dilation theory}, preprint

\bibitem{ACCGM} {A. Avil\'{e}s, F. Cabello S\'{a}nchez, 
J. Castillo, M. Gonz\'{a}lez, Y. Moreno}, {\it Banach
spaces of universal disposition}, J. Funct. Anal. {\bf 261}
(2011), 2347--2361

\bibitem{CsGwK} F. Cabello S\'anchez, J. Garbuli\'nska-W\c{e}grzyn, W. Kubi\'s, {\it Quasi-Banach spaces of almost universal disposition}, preprint

\bibitem{Caradus}
S.R. Caradus, {\it Universal operators and invariant subspaces}, Proc. Amer. Math. Soc. {\bf23} (1969) 526--527

\bibitem{kubis-gar} J. Garbuli\'nska, W. Kubi\'s, {\it Remarks on \Gurarii\ spaces}, Extracta Math. {\bf26} (2011) 235--269

\bibitem{gurarii} V. I. \Gurarii, {\it Space of universal disposition, isotropic spaces and the Mazur problem on rotations of Banach spaces}, Siberian Math. J. {\bf7} (1966), 799--807


\bibitem{KS} W. Kubi\'s, S. Solecki, {\it A proof of uniqueness of the \Gurarii\ space}, Israel J. Math. {\bf195} (2013) 449--456

\bibitem{LinPel}
J. Lindenstrauss, A. Pe\l czy\'nski,
{\it Absolutely summing operators in $L_p$-spaces and their applications},
Studia Math. {\bf29} (1968) 275--326

\bibitem{lusky} W. Lusky, {\it The Gurarij spaces are unique\/}, Arch. Math. (Basel) {\bf27} (1976) 627--635

\bibitem{pelczynski} A. Pe\l{}czy\'nski, {\it Projections in certain Banach spaces}, Studia Math. {\bf19} (1960) 209--228

\end{thebibliography}
\end{document}